# Linearization of the inverse conductivity problem

D V Ingerman

**Abstract.** A positive function (conductivity) on the edges of a graph induces the Dirichlet-to-Neumann map between boundary values of harmonic functions. The inverse conductivity problem is to find the conductivity from the Dirichlet-to-Neumann map. We will show that the map from logarithm of conductivity to the certain logarithms of the determinants of the submatrices of the Dirichlet-to-Neumann map is linear(!) and so the solution of the inverse problem is reduced to solution of the system of linear equations that arise from disjoint paths in the graph. We will make a calculation for a simple tensor product lattice graph and conjecture that it generalizes to planar and three dimensional graphs and also to the continuous case. Depending on the graph the algorithm resembles or not the layer-stripping.

## 1. Dirichlet-to-Neumann map of a graph

The Dirichlet-to-Neumann map is a map from boundary potential to boundary current of a body interior potential satisfying Laplace equation. These maps are the data for inverse boundary problems arising in the problems requiring non-intrusive testing in medicine and oil and gas production industry. There are discrete (on graphs) and continuous (on manifolds) models for the inverse problems.

We will define now the Dirichlet-to-Neumann map of the following graph (network) with 12 edges (conductors) with the positive function (conductivity) $\gamma$ on the edges, 8 boundary vertices (nodes), and 4 interior vertices

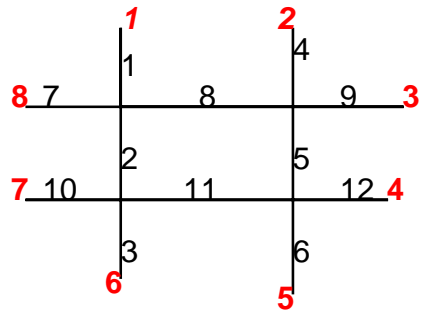

**Figure 1.** A lattice network with 8 boundary nodes, 4 interior nodes and 12 edges.

The Kirchhoff matrix $K$ of the graph is the *12 by 12* matrix with

$$K_{ij} = -\gamma_{ij} \quad \text{if} \quad i \neq j$$

and

$$K_{ii} = \sum_{j} \gamma_{ij}$$

For the lattice network in figure *1*:

$$K = \begin{bmatrix} \gamma_1 & 0 & 0 & 0 & 0 & 0 & 0 & 0 & -\gamma_1 & 0 & & \\ 0 & \gamma_4 & 0 & 0 & 0 & 0 & 0 & 0 & 0 & -\gamma_4 & & \\ 0 & 0 & \gamma_9 & 0 & 0 & 0 & 0 & 0 & 0 & -\gamma_9 & & \\ 0 & 0 & 0 & \gamma_{12} & 0 & 0 & 0 & 0 & 0 & 0 & & \\ 0 & 0 & 0 & 0 & \gamma_6 & 0 & 0 & 0 & 0 & 0 & & \\ 0 & 0 & 0 & 0 & 0 & \gamma_3 & 0 & 0 & 0 & 0 & & \\ 0 & 0 & 0 & 0 & 0 & 0 & \gamma_{10} & 0 & 0 & 0 & & \\ 0 & 0 & 0 & 0 & 0 & 0 & 0 & \gamma_7 & -\gamma_7 & 0 & & \\ -\gamma_1 & 0 & 0 & 0 & 0 & 0 & 0 & -\gamma_7 & \gamma_1+\gamma_7+\gamma_8+\gamma_2 & -\gamma_8 & & \\ 0 & -\gamma_4 & -\gamma_9 & 0 & 0 & 0 & 0 & 0 & -\gamma_8 & \gamma_4+\gamma_9+\gamma_8+\gamma_5 & & \\ 0 & 0 & 0 & -\gamma_{12} & -\gamma_6 & 0 & 0 & 0 & 0 & -\gamma_5 & \gamma_5+\gamma_6 & \\ 0 & 0 & 0 & 0 & 0 & -\gamma_3 & -\gamma_{10} & 0 & -\gamma_2 & 0 & & \end{bmatrix}$$

Numbering boundary vertices of the graph first we write the *12 by 12* Kirchhoff matrix *K* in the block form

$$K = \begin{bmatrix} A & B \\ B^T & C \end{bmatrix}$$

where A is *8 by 8*, B is *8 by 4*, $B^T$ is *4 by 8* and C is *4 by 4*.

A vector *u* on the vertices of the graph is called *harmonic* or *γ-harmonic* if its value at an interior vertex $u_{\text{interior}}$ is the weighted average of its values at the neighbors.

$$u_i = \frac{\sum_j u_j \gamma_{ij}}{\sum_j \gamma_{ij}}$$

or

$$\sum_j \gamma_{ij}(u_i - u_j) = 0$$

(analog of the Laplace-Beltrami equation

$$u_{xx} + u_{yy} = 0$$

or

$$(\gamma u_x)_x + (\gamma u_y)_y = 0$$

.)

The Dirichlet-to-Neumann map *Λ* of the graph is the linear map from the boundary values $u_\partial$ of harmonic functions to the corresponding boundary current.

$$\begin{bmatrix} A & B \\ B^T & C \end{bmatrix} \begin{bmatrix} u_\partial \\ u_{\text{interior}} \end{bmatrix} = \begin{bmatrix} \Lambda u_\partial \\ 0 \end{bmatrix}.$$

Solving we get that the Dirichlet-to-Neumann map *Λ* is *8 by 8* matrix that equals the *Schur complement* of *K* with respect to *C*

$$\Lambda = K/C = A - BC^{-1}B^T.$$

*The inverse conductivity problem* is to find conductivities $\gamma_1, \gamma_2, \ldots, \gamma_{12}$ from the Dirichlet-to-Neumann map $\Lambda$.

## 2. Solution of the inverse conductivity problem

It follows from [1,5,8] that the map

$$\gamma \to \Lambda$$

is invertible and birational. In this paper we reduce the inverse problem to analyzing disjoint paths in the graph and solving a system of linear equations and conjecture that our solution generalizes to planar and three dimensional graphs and also to the continuous case [6].

By matrix block multiplication

$$\begin{bmatrix} A & B \\ B^T & C \end{bmatrix} \begin{bmatrix} I & 0 \\ -C^{-1}B^T & I \end{bmatrix} = \begin{bmatrix} A - BC^{-1}B^T & B \\ 0 & C \end{bmatrix}$$

and therefore

$$\det \begin{bmatrix} A & B \\ B^T & C \end{bmatrix} = \det(A - BC^{-1}B^T) \det C$$

For an $n$ by $n$ matrix $M$ and two subsets $P$ and $Q$ of $\{1, 2, \ldots n\}$ let $M(P,Q)$ denote the submatrix of $M$ consisting of intersection of rows $P$ and columns $Q$ of M.

Let $P$ and $Q$ be two subsets of $\{1,2,3,4,5,6,7,8\}$ and let $I = \{9,10,11,12\}$ then

$$\det \Lambda(P,Q) = \frac{\det K(P \cup I, Q \cup I)}{\det K(I,I)}$$

and expanding the determinants, see [2,5]

$$\det \Lambda(P,Q) = \sum_{disjoint\ pathfromPtoQ} (-1)^{pathsign} \prod_{ij in path} \gamma_{ij} \det K(notinpath, notinpath) / \det K(I,I)$$

For example

$$-\det \Lambda(1,2;5,6) = \gamma_1\gamma_2\gamma_3\gamma_4\gamma_5\gamma_6 / \det K(I,I)$$

$$-\det \Lambda(1,2,8;5,6,8) = \gamma_1\gamma_2\gamma_3\gamma_4\gamma_5\gamma_6\gamma_7 / \det K(I,I)$$

$$\det \Lambda(1;5) = -\frac{\gamma_1\gamma_2\gamma_{11}\gamma_6 \det K(10,10) + \gamma_1\gamma_8\gamma_5\gamma_6 \det K(12,12)}{\det K(I,I)} =$$

$$-\frac{\gamma_1\gamma_2\gamma_{11}\gamma_6(\gamma_8 + \gamma_4 + \gamma_9 + \gamma_5) + \gamma_1\gamma_8\gamma_5\gamma_6(\gamma_2 + \gamma_3 + \gamma_{10} + \gamma_{11})}{\det K(I,I)}$$

$$\det \Lambda(1,5;2,6) = (\gamma_1\gamma_2\gamma_3\gamma_4\gamma_5\gamma_6 - \gamma_1\gamma_8\gamma_4\gamma_3\gamma_{11}\gamma_6)/\det K(I,I)$$

The equations can be represented by the following diagrams if the disjoint paths are unique.

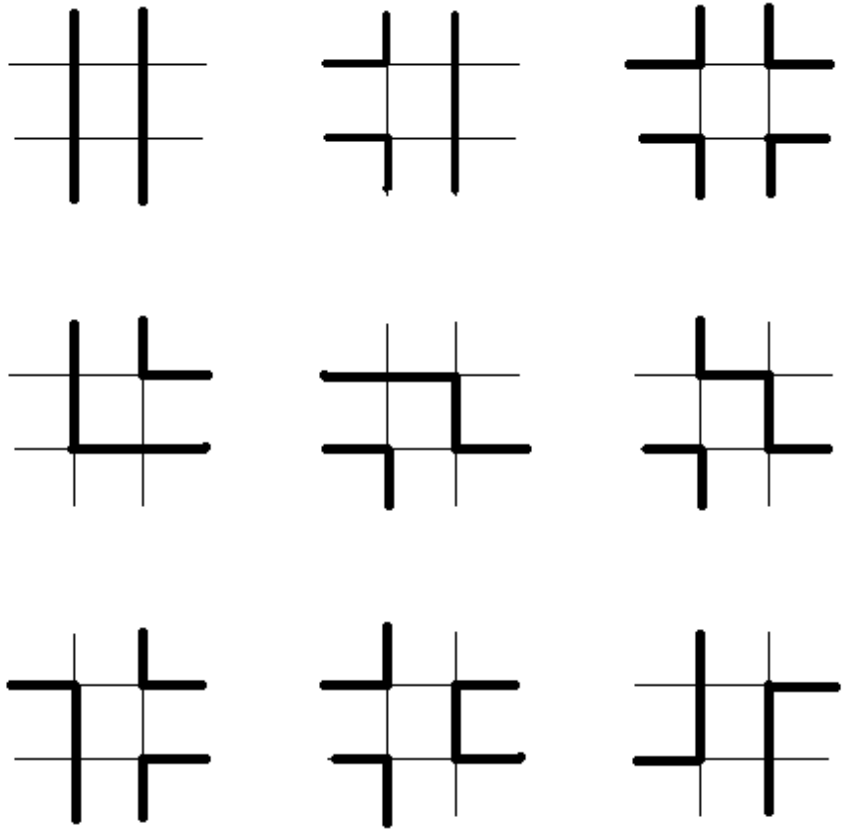

**Figure 2.** Disjoint paths in the network giving linear equations for the inverse problem.

The equations are linear in *log γ* and *log detK(I,I)* and if the disjoint paths go through all four interior vertices of the graph the equations can be written as the following linear system with 13 unknowns

$$\begin{bmatrix} 1 & 1 & 1 & 1 & 1 & 1 & 0 & 0 & 0 & 0 & 0 & 0 & -1 \\ 1 & 1 & 1 & 1 & 1 & 1 & 1 & 0 & 0 & 0 & 0 & 0 & -1 \\ 1 & 0 & 1 & 1 & 1 & 1 & 1 & 0 & 0 & 1 & 0 & 0 & -1 \\ & & & & & & & & & & & & -1 \\ & & & & & & & & & & & & -1 \\ & & & & & & & & & & & & -1 \\ \ldots & \ldots & \ldots & \ldots & \ldots & \ldots & \ldots & \ldots & \ldots & \ldots & \ldots & \ldots & -1 \end{bmatrix} \begin{bmatrix} \log \gamma_1 \\ \log \gamma_2 \\ \log \gamma_3 \\ \log \gamma_4 \\ \log \gamma_5 \\ \log \gamma_6 \\ \log \gamma_7 \\ \log \gamma_8 \\ \log \gamma_9 \\ \log \gamma_{10} \\ \log \gamma_{11} \\ \log \gamma_{12} \\ \log \det K(I,I) \end{bmatrix} = \begin{bmatrix} \log |\det \Lambda(1,2;5,6)| \\ \log |\det \Lambda(1,2,8;5,6,8)| \\ \log |\det \Lambda(1,2,7;8,5,6)| \\ \\ \\ \\ \ldots \end{bmatrix}$$

Most of these equations are linearly dependent but the system is solvable. Subtracting we get with rotations

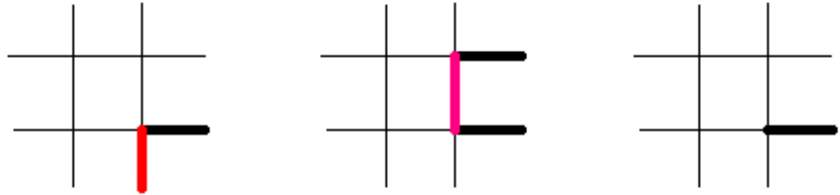

**Figure 3.** Edges in the network corresponding to the non-zero coefficients in the linear system.

where black is + and red is -.

$$\begin{bmatrix} 0 & 0 & 0 & 0 & 0 & 0 & 1 & 0 & 0 & 0 & 0 & 0 & 0 \\ 0 & -1 & 0 & 0 & 0 & 0 & 1 & 0 & 0 & 1 & 0 & 0 & 0 \\ \cdots & \cdots & \cdots & \cdots & \cdots & \cdots & \cdots & \cdots & \cdots & \cdots & \cdots & \cdots & \cdots \\ \cdots & \cdots & \cdots & \cdots & \cdots & \cdots & \cdots & \cdots & \cdots & \cdots & \cdots & \cdots & \cdots \end{bmatrix} \begin{bmatrix} \log \gamma_1 \\ \log \gamma_2 \\ \log \gamma_3 \\ \log \gamma_4 \\ \log \gamma_5 \\ \log \gamma_6 \\ \log \gamma_7 \\ \log \gamma_8 \\ \log \gamma_9 \\ \log \gamma_{10} \\ \log \gamma_{11} \\ \log \gamma_{12} \\ \log \det K(I,I) \end{bmatrix} = \begin{bmatrix} \log|\det \Lambda(1,2,8;5,6,8)| - \log|\det \Lambda(1, \\ \log|\det \Lambda(1,2,7;8,5,6)| - \log|\det \Lambda(1, \\ \cdots \\ \cdots \end{bmatrix}$$

and we easily solve the inverse problem by expressing $\log \gamma$ as a linear combination of $\log|\det \Lambda(P,Q)|$ for various *P* and *Q*.

For an arbitrary graph the rank of the *{1, 0, -1}* linear system arising from the disjoint paths is a fundamental value and requires a special analysis.